\documentclass[a4paper,11pt]{article}
\usepackage{url}
\usepackage[utf8]{inputenc}
\usepackage{amsmath,accents,amssymb,amsthm}
\usepackage[usenames]{color}
\usepackage{longtable}
\usepackage[table]{xcolor}
\usepackage{algpseudocode}
\usepackage{amsfonts}
\usepackage{pgfplots}
\usepackage{subcaption}

\usepackage[normalem]{ulem}
\pgfplotsset{compat=1.18}

\newcommand\mut[1]{\ignorespaces}

\theoremstyle{plain}
\newtheorem{theorem}{Theorem}[section]
\newtheorem{lemma}[theorem]{Lemma}

\newtheorem{corollary}[theorem]{Corollary}
\newtheorem{remark}[theorem]{Remark}
\theoremstyle{definition}

\usepackage[mathscr]{eucal}
\usepackage[usenames]{color}

\usepackage{algorithm,algpseudocode}
\algnewcommand{\Or}{\textbf{ or }}
\algnewcommand{\And}{\textbf{ and }}
\algnewcommand{\From}{\textbf{ from }}
\algnewcommand{\To}{\textbf{ to }}
\algnewcommand{\Let}{\textbf{Let }}
\algnewcommand{\Var}{\texttt}
\algnewcommand{\Int}{\textbf{ int }}
\algnewcommand{\List}{\textbf{ list }}
\algnewcommand{\Bool}{\textbf{ bool }}
\algnewcommand{\True}{\textsc{ true }}
\algnewcommand{\False}{\textsc{ false }}
\renewcommand\leq\leqslant
\renewcommand\geq\geqslant

\usepackage{tikz,tikz-qtree,tikz-qtree-compat,adjustbox,xcolor,array,hhline}
\usepackage{pst-plot}
\usepackage{etoolbox}
\usetikzlibrary{positioning}
\providecommand\circledcolorednumb{}\renewcommand\circledcolorednumb[2]{\resizebox{0.101111\textwidth}{!}{\tikz[baseline=(char.center)]{\node[shape = circle,draw, inner sep = 2pt,fill=#1](char)    {\phantom{00}};\node[anchor=center] at (char.center) {\makebox(0,0){\large{{\sf #2}}}};}}}

\robustify{\circledcolorednumb}

\title{Multiparameter counting of numerical semigroups by Frobenius number, genus, and multiplicity}
\author{Maria Bras-Amorós and Vicenç Torra}

\begin{document}
\maketitle

\begin{abstract}
  A numerical semigroup is a subset of the nonnegative integers, closed under addition and with finite complement. The size of the complement is its genus. The problem of counting semigroups by the largest gap got very important through the so-called Frobenius problem, first documented in 1884. Accordingly, the largest gap is called the Frobenius number. In the last two decades, counting by the genus has become a subject of even more intense study, mostly because of the non-solved conjectures on its monotonic and super-Fibonacci growth.
  It is known that the number of semigroups of each genus grows asymptotically with the genus like the Fibonacci numbers and that the number of semigroups of each Frobenius number grows asymptotically in such a way that each number doubles the previous but one. We prove a formula for the number of numerical semigroups of each Frobenius number, genus, and multiplicity (first nonzero nongap), under the assumption that the multiplicity $m$ and the Frobenius number $F$ satisfy $m\geq\frac{F+1}{3}$. It has been proved that asymptotically almost all semigroups satisfy this condition. From this formula we prove multiparameter recurrence relations. They give, under the required restriction, a multiparameter exact version of the increasing behaviours previously mentioned. Indeed, we prove that the number $n_{F,m,g}$ of semigroups with Frobenius number $F$, multiplicity $m$ and genus $g$ satisfies $n_{F,m,g}=n_{F-2,m-1,g-1}+n_{F-2,m-1,g-2}$ and that the number $n_{F,m}$ of semigroups with Frobenius number $F$ and multiplicity $m$ satisfies $n_{F,m}=2n_{F-2,m-1}.$
\end{abstract}

\begin{center}
{\small {\bf Keywords:}
  Numerical semigroup, genus, Frobenius number, multiplicity
}

{\small {\bf 2020 Mathematics Subject Classification:} 05A15, 06F05, 20M13, 20M14}
\end{center}

\section{Introduction}

A {\em numerical semigroup} is a cofinite monoid of the nonnegative integers. The cardinality of its complement is called the {\em genus}, the first nonzero nongap is called the {\em multiplicity}, and the maximum gap is the {\em Frobenius number} of the numerical semigroup. There have been many efforts to count numerical semigroups by the Frobenius number \cite{fundamentalgaps,DK,DKM} and by the genus \cite{NivaldoMedeiros, Komeda, Br:fibonacci, FromentinHivert, seeds1, rgd, seeds2, DelgadoEliahouFromentin, unleaved}.
The number of semigroups of genus $g$ is usually denoted $n_g$. Zhai \cite{Zhai} proved a conjecture in \cite{Br:fibonacci}, according to which, $n_g$ grows asymptotically like the Fibonacci numbers. On the other side, Backelin \cite{Backelin} proved that the number of semigroups of each Frobenius number grows asymptotically in such a way that each number doubles the previous but one. Let $n_{F,m,g}$ be the number of numerical semigroups with Frobenius number $F$, multiplicity $m$, and genus $g$.
We prove a formula for $n_{F,m,g}$ under the assumption that $m\geq\frac{F+1}{3}$.
It is well known that as the genus grows to infinity, the expected Frobenius number of a typical numerical semigroup approaches twice the multiplicity of the semigroup \cite{KaplanYe,Singhal_Frobenius_genus}. The same behavior is proved when the Frobenius number grows to infinity \cite{Backelin}. Hence, the range where our formula is satisfied covers asymptotically most numerical semigroups, either of a given Frobenius number or of a given genus.
From our formula we derive two multiparameter recurrences which correspond to
the two increasing behaviors previously mentioned. 
Indeed, let $n_{F,m}$ be the number of numerical semigroups 
with Frobenius number $F$ and multiplicity $m$.
On the one hand we prove that,
if $m\geq \frac{F+2}{3}$, then \begin{equation}\label{recurrenceF}n_{F,m}=2n_{F-2,m-1}.\end{equation}
On the other hand,
we prove that,
if $m\geq\frac{F+2}{3}$,
then \begin{equation}\label{recurrenceg}n_{F,m,g}=n_{F-2,m-1,g-1}+n_{F-2,m-1,g-2}.\end{equation} A consequence of the recurrence \eqref{recurrenceg} is the result by Kaplan~\cite{Kaplan} that states that, if $N(m,g)$ is the number of numerical semigroups of multiplicity $m$ and genus $g$, then, under the hypothesis that $m$ and $g$ satisfy $2g \leq 3m-1$, it holds
$N (m, g) = N (m - 1, g - 1) + N (m - 1, g - 2).$ 

The structure of the work is as follows.
In
Section~\ref{s:limits} we stablish the limits outside which $n_{F,m,g}$ is zero and characterize the parameters $F,m,g$ for which $n_{F,m,g}$ is nonzero.
In 
Section~\ref{s:representation} we give a representation of the numerical semigroups whose multiplicity $m$ and Frobenius number $F$ satisfy $m\geq \frac{F+1}{3}$.
In Section~\ref{s:counting}
from the previous representation we deduce the formula for
$n_{F,m,g}$
under the assumption that $m\geq \frac{F+1}{3}$.
In Section~\ref{s:Fgrowth}
we deduce the recurrence \eqref{recurrenceF}, and in
Section~\ref{s:ggrowth}
we deduce the recurrence \eqref{recurrenceg}.

\section{Multiparameter counting by the Frobenius number, multiplicity, and genus}\label{s:limits}

Let ${\mathcal S}_{F,m}$ be the set of numerical semigroups of Frobenius number $F$ and multiplicity $m$ and let ${\mathcal S}_{F,m,g}$ be the set of numerical semigroups of Frobenius number $F$, multiplicity $m$, and genus $g$.
We are interested in the cardinalities of the sets ${\mathcal S}_{F,m}$ and ${\mathcal S}_{F,m,g}$, which we called $n_{F,m}$ and $n_{F,m,g}$, respectively.
As an illustrative example, in Table~\ref{tF24} we give the table with all the nonzero values $n_{F=24,m,g}$, together with a right column with all the nonzero values $n_{F=24,m}$. Similarly, in Table~\ref{tF25} and Table~\ref{tF26} we give the corresponding tables with Frobenius number 25 and 26. They will be used as paradigmatic examples throughout the paper.

\begin{table}
{\centering
  \resizebox{.7\textwidth}{!}{\begin{tabular}{||c||c|c|c|c|c|c|c|c|c|c|c|c||c||}
\hline\hline
 & $g=$13 & $g=$14 & $g=$15 & $g=$16 & $g=$17 & $g=$18 & $g=$19 & $g=$20 & $g=$21 & $g=$22 & $g=$23 & $g=$24 & total\\\hline\hline
$m=2$ &  &  &  &  &  &  &  &  &  &  &  &  & $0$ \\\hline
$m=3$ &  &  &  &  &  &  &  &  &  &  &  &  & $0$ \\\hline
$m=4$ &  &  &  &  &  &  &  &  &  &  &  &  & $0$ \\\hline
$m=5$ & $2$ & $5$ & $8$ & $9$ & $8$ & $6$ & $3$ & $1$ &  &  &  &  & $42$ \\\hline
$m=6$ &  &  &  &  &  &  &  &  &  &  &  &  & $0$ \\\hline
$m=7$ & $4$ & $12$ & $23$ & $32$ & $31$ & $23$ & $13$ & $5$ & $1$ &  &  &  & $144$ \\\hline
$m=8$ &  &  &  &  &  &  &  &  &  &  &  &  & $0$ \\\hline
$m=9$ & $4$ & $19$ & $45$ & $73$ & $87$ & $76$ & $49$ & $23$ & $7$ & $1$ &  &  & $384$ \\\hline
$m=10$ & $2$ & $14$ & $44$ & $83$ & $106$ & $97$ & $64$ & $29$ & $8$ & $1$ &  &  & $448$ \\\hline
$m=11$ & $1$ & $9$ & $36$ & $84$ & $126$ & $126$ & $84$ & $36$ & $9$ & $1$ &  &  & $512$ \\\hline
$m=12$ &  &  &  &  &  &  &  &  &  &  &  &  & $0$ \\\hline
$m=13$ & $1$ & $10$ & $45$ & $120$ & $210$ & $252$ & $210$ & $120$ & $45$ & $10$ & $1$ &  & $1024$ \\\hline
$m=14$ &  & $1$ & $9$ & $36$ & $84$ & $126$ & $126$ & $84$ & $36$ & $9$ & $1$ &  & $512$ \\\hline
$m=15$ &  &  & $1$ & $8$ & $28$ & $56$ & $70$ & $56$ & $28$ & $8$ & $1$ &  & $256$ \\\hline
$m=16$ &  &  &  & $1$ & $7$ & $21$ & $35$ & $35$ & $21$ & $7$ & $1$ &  & $128
$ \\\hline
$m=17$ &  &  &  &  & $1$ & $6$ & $15$ & $20$ & $15$ & $6$ & $1$ &  & $64$ \\\hline
$m=18$ &  &  &  &  &  & $1$ & $5$ & $10$ & $10$ & $5$ & $1$ &  & $32$ \\\hline
$m=19$ &  &  &  &  &  &  & $1$ & $4$ & $6$ & $4$ & $1$ &  & $16$ \\\hline
$m=20$ &  &  &  &  &  &  &  & $1$ & $3$ & $3$ & $1$ &  & $8$ \\\hline
$m=21$ &  &  &  &  &  &  &  &  & $1$ & $2$ & $1$ &  & $4$ \\\hline
$m=22$ &  &  &  &  &  &  &  &  &  & $1$ & $1$ &  & $2$ \\\hline
$m=23$ &  &  &  &  &  &  &  &  &  &  & $1$ &  & $1$ \\\hline
$m=24$ &  &  &  &  &  &  &  &  &  &  &  &  & $0$ \\\hline
$m=25$ &  &  &  &  &  &  &  &  &  &  &  & $1$ & $1$ \\\hline
\hline total &
$14$ & $70$ & $211$ & $446$ & $688$ & $790$ & $675$ & $424$ & $190$ & $58$ & $11$ & $1$ & $\Sigma=3578$
\\\hline\hline
\end{tabular}}
    \caption{Values of $n_{F=24,m,g}$ and $n_{F=24,m}$.}
\label{tF24}
}
\end{table}
\begin{table}
{\centering
  \resizebox{.8\textwidth}{!}{\begin{tabular}{||c||c|c|c|c|c|c|c|c|c|c|c|c|c||c||}
\hline\hline
 & $g=$13 & $g=$14 & $g=$15 & $g=$16 & $g=$17 & $g=$18 & $g=$19 & $g=$20 & $g=$21 & $g=$22 & $g=$23 & $g=$24 & $g=$25 & total\\\hline\hline
$m=2$ & $1$ &  &  &  &  &  &  &  &  &  &  &  &  & $1$ \\\hline
$m=3$ & $1$ & $1$ & $1$ & $1$ & $1$ &  &  &  &  &  &  &  &  & $5$ \\\hline
$m=4$ & $4$ & $5$ & $4$ & $4$ & $3$ & $2$ & $1$ &  &  &  &  &  &  & $23$ \\\hline
$m=5$ &  &  &  &  &  &  &  &  &  &  &  &  &  & $0$ \\\hline
$m=6$ & $8$ & $16$ & $25$ & $27$ & $24$ & $17$ & $10$ & $4$ & $1$ &  &  &  &  & $132$ \\\hline
$m=7$ & $6$ & $16$ & $31$ & $42$ & $46$ & $40$ & $27$ & $14$ & $5$ & $1$ &  &  &  & $228$ \\\hline
$m=8$ & $8$ & $26$ & $54$ & $76$ & $82$ & $68$ & $42$ & $19$ & $6$ & $1$ &  &  &  & $382$ \\\hline
$m=9$ & $8$ & $31$ & $71$ & $115$ & $140$ & $134$ & $99$ & $56$ & $24$ & $7$ & $1$ &  &  & $686$ \\\hline
$m=10$ & $4$ & $24$ & $69$ & $128$ & $170$ & $168$ & $126$ & $72$ & $30$ & $8$ & $1$ &  &  & $800$ \\\hline
$m=11$ & $2$ & $16$ & $58$ & $127$ & $189$ & $203$ & $161$ & $93$ & $37$ & $9$ & $1$ &  &  & $896$ \\\hline
$m=12$ & $1$ & $10$ & $45$ & $120$ & $210$ & $252$ & $210$ & $120$ & $45$ & $10$ & $1$ &  &  & $1024$ \\\hline
$m=13$ & $1$ & $11$ & $55$ & $165$ & $330$ & $462$ & $462$ & $330$ & $165$ & $55$ & $11$ & $1$ &  & $2048$ \\\hline
$m=14$ &  & $1$ & $10$ & $45$ & $120$ & $210$ & $252$ & $210$ & $120$ & $45$ & $10$ & $1$ &  & $1024$ \\\hline
$m=15$ &  &  & $1$ & $9$ & $36$ & $84$ & $126$ & $126$ & $84$ & $36$ & $9$ & $1$ &  & $512$ \\\hline
$m=16$ &  &  &  & $1$ & $8$ & $28$ & $56$ & $70$ & $56$ & $28$ & $8$ & $1$ &  & $256$ \\\hline
$m=17$ &  &  &  &  & $1$ & $7$ & $21$ & $35$ & $35$ & $21$ & $7$ & $1$ &  & $128$ \\\hline
$m=18$ &  &  &  &  &  & $1$ & $6$ & $15$ & $20$ & $15$ & $6$ & $1$ &  & $64$ \\\hline
$m=19$ &  &  &  &  &  &  & $1$ & $5$ & $10$ & $10$ & $5$ & $1$ &  & $32$ \\\hline
$m=20$ &  &  &  &  &  &  &  & $1$ & $4$ & $6$ & $4$ & $1$ &  & $16$ \\\hline
$m=21$ &  &  &  &  &  &  &  &  & $1$ & $3$ & $3$ & $1$ &  & $8$ \\\hline
$m=22$ &  &  &  &  &  &  &  &  &  & $1$ & $2$ & $1$ &  & $4$ \\\hline
$m=23$ &  &  &  &  &  &  &  &  &  &  & $1$ & $1$ &  & $2$ \\\hline
$m=24$ &  &  &  &  &  &  &  &  &  &  &  & $1$ &  & $1$ \\\hline
$m=25$ &  &  &  &  &  &  &  &  &  &  &  &  &  & $0$ \\\hline
$m=26$ &  &  &  &  &  &  &  &  &  &  &  &  & $1$ & $1$ \\\hline
\hline total &
$44$ & $157$ & $424$ & $860$ & $1360$ & $1676$ & $1600$ & $1170$ & $643$ & $256$ & $70$ & $12$ & $1$ & $\Sigma=8273$
\\\hline\hline
\end{tabular}}
  \caption{Values of $n_{F=25,m,g}$ and $n_{F=25,m}$.}
\label{tF25}
}
\end{table}
\begin{table}
{\centering
  \resizebox{.9\textwidth}{!}{\begin{tabular}{||c||c|c|c|c|c|c|c|c|c|c|c|c|c||c||}
\hline\hline
 & $g=$14 & $g=$15 & $g=$16 & $g=$17 & $g=$18 & $g=$19 & $g=$20 & $g=$21 & $g=$22 & $g=$23 & $g=$24 & $g=$25 & $g=$26 & total\\\hline\hline
$m=2$ &  &  &  &  &  &  &  &  &  &  &  &  &  & $0$ \\\hline
$m=3$ & $1$ & $1$ & $1$ & $1$ & $1$ &  &  &  &  &  &  &  &  & $5$ \\\hline
$m=4$ & $1$ & $2$ & $3$ & $4$ & $3$ & $2$ & $1$ &  &  &  &  &  &  & $16$ \\\hline
$m=5$ & $3$ & $6$ & $8$ & $9$ & $8$ & $6$ & $3$ & $1$ &  &  &  &  &  & $44$ \\\hline
$m=6$ & $3$ & $10$ & $18$ & $23$ & $23$ & $17$ & $10$ & $4$ & $1$ &  &  &  &  & $109$ \\\hline
$m=7$ & $7$ & $17$ & $28$ & $39$ & $42$ & $37$ & $26$ & $14$ & $5$ & $1$ &  &  &  & $216$ \\\hline
$m=8$ & $4$ & $18$ & $42$ & $65$ & $74$ & $64$ & $41$ & $19$ & $6$ & $1$ &  &  &  & $334$ \\\hline
$m=9$ & $8$ & $28$ & $66$ & $109$ & $136$ & $133$ & $99$ & $56$ & $24$ & $7$ & $1$ &  &  & $667$ \\\hline
$m=10$ & $4$ & $23$ & $64$ & $118$ & $160$ & $163$ & $125$ & $72$ & $30$ & $8$ & $1$ &  &  & $768$ \\\hline
$m=11$ & $2$ & $16$ & $58$ & $127$ & $189$ & $203$ & $161$ & $93$ & $37$ & $9$ & $1$ &  &  & $896$ \\\hline
$m=12$ & $1$ & $10$ & $45$ & $120$ & $210$ & $252$ & $210$ & $120$ & $45$ & $10$ & $1$ &  &  & $1024$ \\\hline
$m=13$ &  &  &  &  &  &  &  &  &  &  &  &  &  & $0$ \\\hline
$m=14$ & $1$ & $11$ & $55$ & $165$ & $330$ & $462$ & $462$ & $330$ & $165$ & $55$ & $11$ & $1$ &  & $2048$ \\\hline
$m=15$ &  & $1$ & $10$ & $45$ & $120$ & $210$ & $252$ & $210$ & $120$ & $45$ & $10$ & $1$ &  & $1024$ \\\hline
$m=16$ &  &  & $1$ & $9$ & $36$ & $84$ & $126$ & $126$ & $84$ & $36$ & $9$ & $1$ &  & $512$ \\\hline
$m=17$ &  &  &  & $1$ & $8$ & $28$ & $56$ & $70$ & $56$ & $28$ & $8$ & $1$ &  & $256$ \\\hline
$m=18$ &  &  &  &  & $1$ & $7$ & $21$ & $35$ & $35$ & $21$ & $7$ & $1$ &  & $128$ \\\hline
$m=19$ &  &  &  &  &  & $1$ & $6$ & $15$ & $20$ & $15$ & $6$ & $1$ &  & $64$ \\\hline
$m=20$ &  &  &  &  &  &  & $1$ & $5$ & $10$ & $10$ & $5$ & $1$ &  & $32$ \\\hline
$m=21$ &  &  &  &  &  &  &  & $1$ & $4$ & $6$ & $4$ & $1$ &  & $16$ \\\hline
$m=22$ &  &  &  &  &  &  &  &  & $1$ & $3$ & $3$ & $1$ &  & $8$ \\\hline
$m=23$ &  &  &  &  &  &  &  &  &  & $1$ & $2$ & $1$ &  & $4$ \\\hline
$m=24$ &  &  &  &  &  &  &  &  &  &  & $1$ & $1$ &  & $2$ \\\hline
$m=25$ &  &  &  &  &  &  &  &  &  &  &  & $1$ &  & $1$ \\\hline
$m=26$ &  &  &  &  &  &  &  &  &  &  &  &  &  & $0$ \\\hline
$m=27$ &  &  &  &  &  &  &  &  &  &  &  &  & $1$ & $1$ \\\hline
\hline total &
$35$ & $143$ & $399$ & $835$ & $1341$ & $1669$ & $1600$ & $1171$ & $643$ & $256$ & $70$ & $12$ & $1$ & $\Sigma=8175$
\\\hline\hline
\end{tabular}}
  \caption{Values of $n_{F=26,m,g}$ and $n_{F=26,m}$.}
\label{tF26}
}
\end{table}

Let $[a,b]=\{i\in{\mathbb Z}\mbox{ such that }a\leq i\leq b\}$ and
$[a,\infty)=\{i\in{\mathbb Z}\mbox{ such that }i\geq a\}$.
The next lemma analyzes the semigroups for which all elements smaller than the Frobenius number are multiples of the multiplicity. The result will be used then in Lemma~\ref{polygon}.

\begin{lemma}\label{spine}
  Suppose $m$ is an integer, $m\geq 2$, that does not divide $F$.
  Consider the set $$S_r=\{0,m,2m,\dots,r\cdot m\}\cup[F+1,\infty).$$
    \begin{enumerate}
    \item The set $S_r$ is a numerical semigroup if and only if $r\geq \lfloor\frac{F}{m}\rfloor$.
    \item The set $S_r$ with $r=\lfloor\frac{F}{m}\rfloor$ is the unique element in ${\mathcal S}_{F,m,(F-\lfloor\frac{F}{m}\rfloor)}$.
    \end{enumerate}
\end{lemma}

\begin{proof}
  \begin{enumerate}
    \item First we notice that $r\geq \lfloor\frac{F}{m}\rfloor$ is equivalent to $r\cdot m+m\geq F+1$.
If $S_r$ is a numerical semigroup, then, necessarily, since $m, r\cdot m\in S_r$, it must be $r\cdot m+m\in S_r$, and, since $r\cdot m+m>r\cdot m$, then $r\cdot m+m\geq F+1$. 
Conversely, if $r\cdot m+m\geq F+1$, then all multiples of $m$ belong to $S_r$ and, since all elements smaller than $F$ are multiples of $m$, the set $S_r$ is a numerical semigroup.
\item If $r=\lfloor\frac{F}{m}\rfloor$ then, by the previous statement, $S_r$ is a numerical semigroup. Its multiplicity and Frobenius number are obviously $m$ and $F$, respectively. The genus of $S_r$ is $F-\lfloor\frac{F}{m}\rfloor$ because $\{0,m,2m,\dots,r\cdot m\}\cap[F+1,\infty)=\emptyset$. Hence, 
  $S_r\in {\mathcal S}_{F,m,(    F-\lfloor\frac{F}{m}\rfloor    )}$.

  Now, suppose that a semigroup $S$ belongs to ${\mathcal S}_{F,m,(F-\lfloor\frac{F}{m}\rfloor)}$.
  Since $m\in S$, we have $\{0,m,2m,\dots,r\cdot m\}\subseteq S$. Since the Frobenius number of $S$ is $F$, we have $[F+1,\infty)\subseteq S$. Therefore,
    $S_r\subseteq S$. Finally, since the genus of both $S_r$ and $S$ is the same, necessarily $S=S_r$.
  \end{enumerate}
\end{proof}

The next lemma analyzes the semigroups for which all elements smaller than the Frobenius number are consecutive from $m$ to $F-1$. The result will be used then in Lemma~\ref{polygon}.

\begin{lemma}\label{tableau}
  Suppose that $F,m$ are integers with $2\leq m\leq F-1$.
  Consider the set $$S=\{0\}\cup[m,F-1]\cup[F+1,\infty).$$
    \begin{enumerate}
    \item The set $S$ is a numerical semigroup if and only if $m\geq \frac{F+1}{2}$.
    \item If $m\geq \frac{F+1}{2}$, the set $S$
      is the unique element in ${\mathcal S}_{F,m,m}$.
    \end{enumerate}
\end{lemma}

\begin{proof}
  \begin{enumerate}
  \item
The set $S$ is a numerical semigroup if and only if the sum of two elements in $[m,F-1]$ is in $S$. The set of sums of two elements in $[m,F-1]$ is $[2m,2F-2]$. Since this is an interval with its maximum being larger than $F$, it belongs to $S$ if and only if its minimum is larger than $F$. That is, if and only if $2m\geq F+1$, and the result follows. 
\item If $m\geq \frac{F+1}{2}$, by the previous statement, $S$ is a numerical semigroup. Its multiplicity and Frobenius number are obviously $m$ and $F$, respectively. The genus of $S$ is $F-(F-m)=m$.
  Hence, 
$S\in {\mathcal S}_{F,m,m}$.

  Now, suppose that a semigroup $\tilde S$ belongs to ${\mathcal S}_{F,m,m}$.
  Since $F\not\in \tilde S$ and $m$ is the smallest nonzero element in $\tilde S$, we have $\tilde S\subseteq\{0\}\cup[m,F-1]\cup[F+1,\infty)=S$.
    Now, since the genus of both $\tilde S$ and $S$ is the same, necessarily $\tilde S=S$.
  \end{enumerate}
  \end{proof}

One first result that one can observe looking at the tables is the following lemma.

\begin{lemma}\label{polygon}
  \begin{enumerate}
    \item
      \begin{enumerate}
    \item If $m=1$, then $n_{F,m,g}=0$ for all $F\geq0$ and all $g\geq 0$.
    \item If $g<\frac{F+1}{2}$, then $n_{F,m,g}=0$ for all $m\geq 2$.
    \item If $g=F$, then $n_{F,m,g}=\left\{\begin{array}{ll}1&\mbox{if }m=F+1\\0&\mbox{otherwise}\end{array}\right.$
    \item If $g>F$, then $n_{F,m,g}=0$ for all $m\geq 2$.
    \item If $m=F-1$, then $n_{F,m,g}=\left\{\begin{array}{ll}1&\mbox{if }g=F-1\\0&\mbox{otherwise}\end{array}\right.$
    \item If $m=F$, then $n_{F,m,g}=0$ for all $g\geq 0$.
    \item If $m=F+1$, then $n_{F,m,g}=\left\{\begin{array}{ll}1&\mbox{if }g=F\\0&\mbox{otherwise}\end{array}\right.$
    \item If $m>F+1$, then $n_{F,m,g}=0$ for all $g\geq 0$.
      \end{enumerate}
    \item If $\frac{F+1}{2}\leq g\leq F-1$ and $2\leq m\leq F-1$, then the following hold.
    \begin{enumerate}
    \item If $g=F-\lfloor\frac{F}{m}\rfloor$, then $n_{F,m,g}=1$.
    \item If $g>F-\lfloor\frac{F}{m}\rfloor$, then $n_{F,m,g}=0$.
    \item If $g=m$, then $n_{F,m,g}=1$.
    \item If $g<m$, then $n_{F,m,g}=0$.
    \end{enumerate}
    \end{enumerate}
\end{lemma}

\begin{proof}
  \begin{enumerate} \item \begin{enumerate}\item $m=1$ implies that the Frobenius number is not defined or defined as $-1$.
     \item
       It is a consequence of the fact that any numerical semigroup satisfies $F+1\leq 2g$.
       
     \item 
       It is a consequence of the fact that $\{0\}\cup [F+1,\infty)$ is the unique semigroup with Frobenius number $F$ equal to its genus.
       
    \item It is a consequence of the fact that any semigroup satisfies
      $g\leq F$.
    \item The unique semigroup with $m=F-1$ is $\{0,F-1\}\cup[F+1,\infty)$, which has genus $F-1$.
    \item The multiplicity, which is a nongap, can not be the Frobenius number, which is a gap. Hence $n_{a,a,g}=0$ for any $a$ and any $g$.
    \item The unique semigroup with $m=F+1$ is $\{0\}\cup[F+1,\infty)$, which has genus $F$.
    \item By definition of the Frobenius number, $F+1$ is a nongap, and, so, $F+1\geq m$.
  \end{enumerate}
  \item
Statements $(a)$ and $(b)$ are a consequence of Lemma~\ref{spine}, while statements $(c)$ and $(d)$ are a consequence of Lemma~\ref{tableau}.
  \end{enumerate}
      
    \end{proof}

As a consequence of Lemma~\ref{polygon}, for a fixed Frobenius number $F$, we have that $n_{F,m,g}$ may only be nonzero either for $g=F$ and $m=F+1$, in which case $n_{F,F+1,F}=1$, or inside the polygon delimited by the inequalities (i) $\frac{F+1}{2}\leq g\leq F-1$, (ii) $2\leq m\leq F-1$, and (iii) $m\leq g\leq F-\lfloor\frac{F}{m}\rfloor$.
The next lemma characterizes the tuples $F,m,g$ for which $n_{F,m,g}=0$ inside this polygon.

For general sets of integers $A,B$ and an integer $n$ define $n+A=\{n+a:a\in A\}$, $n-A=\{n-a:a\in A\}$, $A+B=\{a+b:a\in A,b\in B\}$.

\begin{lemma}\label{zerorows}
  Let $F,m,g$ be integers such that $\frac{F+1}{2}\leq g\leq F-1$, $2\leq m\leq F-1$, and $m\leq g\leq F-\lfloor\frac{F}{m}\rfloor$.
  Then $n_{F,m,g}=0$ if and only if $m$ divides $F$.
\end{lemma}

\begin{proof}
  It is obvious that if $m$ divides $F$, then $n_{F,m,g}=0$ for any $g$.
  Now, suppose that $m$ does not divide $F$.
    By hypothesis, $g$ is bounded by $\max\{\frac{F+1}{2},m\}\leq g\leq \min\{F-1, F-\lfloor\frac{F}{m}\rfloor\}$. We will divide the proof into two parts. First, we will consider the case in which $m\geq \frac{F+1}{2}$, in which case $\max\{\frac{F+1}{2},m\}=m$ and $\min\{F-1, F-\lfloor\frac{F}{m}\rfloor\}=F-1$.
    Then, we will consider the case in which $m< \frac{F+1}{2}$, in which case $\max\{\frac{F+1}{2},m\}=\frac{F+1}{2}$ and $\min\{F-1,F-\lfloor\frac{F}{m}\rfloor\}=F-\lfloor\frac{F}{m}\rfloor$.
  
  If $m\geq\frac{F+1}{2}$, then the set $$I_1=\{0,m\}\cup[F+1,\infty)$$ is a numerical semigroup of genus $F-1$ and, by Lemma~\ref{tableau}, the set $$I_2=\{0\}\cup[m,F-1]\cup[F+1,\infty)$$ is a numerical semigroup of genus $m$.
      Any set $S$ such that $I_1\subseteq S\subseteq I_2$ will be a numerical semigroup with multiplicity $m$ and Frobenius number $F$. And there will be at least one such subset for each genus from $m$ to $F-1$. Hence, for $m\geq \frac{F+1}{2}$, we have $n_{F,m,g}\neq 0$ within the limits $m\leq g\leq F-1$.

  If $m<\frac{F+1}{2}$, then, by Lemma~\ref{spine},
  $$J_1=\left\{0,m,2m,\dots,\left\lfloor\frac{F}{m}\right\rfloor m\right\}\cup[F+1,\infty)$$ is a numerical semigroup of genus $F-\lfloor\frac{F}{m}\rfloor$.
    Now define define $J_0=\left\{m,2m,\dots,    \left\lfloor    \frac{\lfloor\frac{F}{2}\rfloor}{m}    \right\rfloor m    \right\}$. We claim that
    \begin{enumerate}
\item[(i)] $F-J_0\subseteq \left[\left\lfloor\frac{F}{2}\right\rfloor+1,F-1\right];$
\item[(ii)] The set $$A=\left[\left\lfloor\frac{F}{2}\right\rfloor+1,F-1\right]\setminus(F-J_0)$$ contains all multiples of $m$ lying in the interval
  $\left[\left\lfloor\frac{F}{2}\right\rfloor+1,F-1\right]$;
  \item[(iii)] The set $A$ does not contain any element congruent to $F$ modulo $m$.
      \end{enumerate}
    Indeed, to see (i), suppose that $j\in J_0$. Then,
    $m\leq j\leq\left\lfloor\frac{\left\lfloor\frac{F}{2}\right\rfloor}{m}\right\rfloor m$. Notice that, if $F$ is odd, then
    $\left\lfloor\frac{\left\lfloor\frac{F}{2}\right\rfloor}{m}\right\rfloor m\leq\lfloor\frac{F}{2}\rfloor=\lfloor\frac{F-1}{2}\rfloor$. On the other hand, if $F$ is even, since $m$ does not divide $F$, then $m$ does not divide $\lfloor\frac{F}{2}\rfloor$ either and, hence, we have
    $\left\lfloor\frac{\left\lfloor\frac{F}{2}\right\rfloor}{m}\right\rfloor m\leq\lfloor\frac{F}{2}\rfloor-1=\lfloor\frac{F-1}{2}\rfloor$.
Hence, if $j\in J_0$, then,
$m\leq j\leq\left\lfloor\frac{F-1}{2}\right\rfloor$.
Now, $F-j$ must be bounded by $F-\left\lfloor\frac{F-1}{2}\right\rfloor\leq F-j\leq F-m$. Since, $F-\left\lfloor\frac{F-1}{2}\right\rfloor=\left\lfloor\frac{F}{2}\right\rfloor+1$, it follows that $F-j\in \left[\left\lfloor\frac{F}{2}\right\rfloor+1,F-1\right]$.

Statement (ii) follows from the fact that $F-J_0$ contains only elements that are congruent to $F$ modulo $m$, while $F$ is not a multiple of $m$ by hypothesis.

To see statement (iii), notice that the elements in
$\left[\left\lfloor\frac{F}{2}\right\rfloor+1,F-1\right]$ congruent to $F$ modulo $m$ are exactly $F-m,F-2m,\dots,F-rm$, with $r$ the maximum integer such that
$F-rm\geq \lfloor\frac{F}{2}\rfloor+1$, that is,
$r=\lfloor\frac{F-1-\lfloor\frac{F}{2}\rfloor}{m}\rfloor=\lfloor\frac{\lceil\frac{F}{2}\rceil-1}{m}\rfloor=\lfloor\frac{\lfloor\frac{F-1}{2}\rfloor}{m}\rfloor.$
If $F$ is odd, it is easy to check that
$\lfloor\frac{\lfloor\frac{F-1}{2}\rfloor}{m}\rfloor=\lfloor\frac{\lfloor\frac{F}{2}\rfloor}{m}\rfloor$, while, if $F$ is even,
$\lfloor\frac{\lfloor\frac{F-1}{2}\rfloor}{m}\rfloor=\lfloor\frac{\frac{F}{2}-1}{m}\rfloor$, which, since $m$ does not divide $\frac{F}{2}$, must equal
$\lfloor\frac{\frac{F}{2}}{m}\rfloor=\lfloor\frac{\lfloor\frac{F}{2}\rfloor}{m}\rfloor.$
Hence, by the definition of $J_0$, the elements in
$\left[\left\lfloor\frac{F}{2}\right\rfloor+1,F-1\right]$ congruent to $F$ modulo $m$ are exactly $F-J_0$, and the result in (iii) follows.

Now, from the previous claim we can derive that the set
    $$J_2=\{0\}\cup J_0\cup A   \cup[F+1,\infty)$$
      is a numerical semigroup.
              Indeed, suppose $a,b\in J_2$, and we need to see that $a+b\in J_2$. If $a$ or $b$ are larger than $F$ then this is obvious. Otherwise, if $a,b\in J_0$, then $a+b$ is a multiple of $m$ and so, it either belongs to $J_0$, or to $A$ (by claim (ii)) or is larger than $F$. Hence, $a+b\in J_2$.
        If $a,b\in A$, then $a+b>F$ and so, $a+b\in J_2$.
        Now suppose that $a\in J_0 $ and $b\in A$. Then $a+b\geq \lfloor\frac{F}{2}\rfloor +1$ and so, we need to see that $a+b$ is either in $A$ or it is larger than $F$. Suppose $a+b\leq F$. By definition of $A$, the number $a+b$ can not be equal to $F$. Hence, $a+b\in[\lfloor\frac{F}{2}\rfloor+1,F-1]$ and so, if $a+b\in J_2$, then $a+b$ must belong to
        $J_2\cap [\lfloor\frac{F}{2}\rfloor+1,F-1]=A$. To see that in this case $a+b\in A$, we need to see that $a+b\not\in F- J_0$. Indeed, if $a+b=F-a'$ with $a'\in J_0$, then $b=F-(a+a')$, which is congruent to $F$ modulo $m$ and so, by claim (iii), $b\not\in A$, a contradiction.

      The sets $\{0\}$, $J_0$, $A$, and $[F+1,\infty)$ are disjoint and so, the genus of the semigroup $J_2$ is $F-\#J_0-\#A$.
          Now, by the definition of the set $A$ and by claim (i), $\#A=F-\lfloor\frac{F}{2}\rfloor-1-\#J_0$ and, consequently, the genus of $J_2$ is $\lfloor\frac{F}{2}\rfloor+1=\lceil\frac{F+1}{2}\rceil$.
            
      Any set $S$ such that $J_1\subseteq S\subseteq J_2$ will be a numerical semigroup. The proof is analogous to the proof of the fact that $J_2$ is a numerical semigroup, just taking into account that all multiples of $m$ belong to $J_1$ and so, to $S$.

      Any semigroup $S$ with $J_1\subseteq S\subseteq J_2$ has multiplicity $m$ and Frobenius number $F$. And there will be at least one such semigroup for each genus from $ F-\lfloor\frac{F}{m}\rfloor$
      to $\lceil\frac{F+1}{2}\rceil$. Hence, for $m< \frac{F+1}{2}$, we have $n_{F,m,g}\neq 0$ within the limits $\frac{F+1}{2}\leq g\leq F-\lfloor\frac{F}{m}\rfloor$.
\end{proof}

Lemma~\ref{zerorows} can be checked in Table~\ref{tF24}, where the rows corresponding to $m=2,3,4,6,8,12$ are empty, while all the other rows are non-empty exactly inside the polygon given by the bounds on $g$.

The conclusion of this section is the next theorem, which is a consequence of Lemma~\ref{polygon} and Lemma~\ref{zerorows}.

\begin{theorem}
  $n_{F_m,g}\neq 0$ if and only if either
  \begin{itemize}
  \item $g=F$ and $m=F+1$,
  \item (i) $\frac{F+1}{2}\leq g\leq F-1$, (ii) $2\leq m\leq F-1$, (iii) $m\leq g\leq F-\lfloor\frac{F}{m}\rfloor$, and (iv) $m$ does not divide $F$.
   \end{itemize}
  \end{theorem}

\section{A representation of the numerical semigroups with $m\geq\frac{F+1}{3}$}
\label{s:representation}

In this section we will study semigroups for which the multiplicity $m$ and the Frobenius number $F$ satisfy $m\geq \frac{F+1}{3}$. We will base all our results on the definition of the sets of subsets $B_d$ that follows. 
We want to remark that this definition is equivalent to that of the set of subsets denoted ${\mathcal A}_k$ in Zhao's reference \cite{Zhao}, with $k=d$. In that reference, Zhao defined a semigroup $\Lambda$ to be of type $(A,k)$, for $A\in{\mathcal A}_k$, if $\Lambda\cap[m,m+k]=m+A$, and showed that every numerical semigroup with $2m<F<3m$ is of type $(A,k)$ for a unique $k$ and $A\in{\mathcal A}_k$. Then he described the set of all semigroups with type $(A,k)$ for $A \in {\mathcal A}_k$, such that $2m<F<3m$, and used it to derive the number of such semigroups with a given genus $g$.
  In a similar but slightly different approach, reference \cite{BrasKaplanSinghal} focuses on the classification of semigroups with $2m<F<3m$ in terms of
  their intersection with $[m,F-m]$ and with $[2m,F]$. In this case, the first intersection would correspond to the set $m+A$ for some $A\in{\mathcal A}_{F-2m}$ according to Zhao's definition.
  We will push Zhao's construction further in order to describe and count all semigroups with a given Frobenius number $F$ and a given multiplicity $m$,
and to count all semigroups with a given combination of parameters $F,m,g$,
whenever $2m\leq F<3m$ (see Corrollary~\ref{cor} and Theorem~\ref{pot2}). These formulas will allow us to derive the multiparameter versions of
the two-step doubling growth of the number of semigroups by Frobenius number
(Theorem~\ref{t:twostepdoubling}) and the Fibonacci-like growth of the number of semigroups by genus (Theorem~\ref{kaplanrefinat}) in the following sections.

For a positive integer $d$, define $B_d$ to be the set of subsets
\begin{eqnarray*}B_d&=&\{A\subseteq[1,d-1]\mbox{ such that }d\not\in A+A\}\\
  &=&\left\{A\subseteq[1,d-1] \mbox{ such that }\{t,d-t\}\not\subseteq A \mbox{ for all }t\in\left[1,\left\lfloor\frac{d}{2}\right\rfloor\right]\right\}
\end{eqnarray*}
      The subtle difference between this definition of $B_d$ and the definition of ${\mathcal A}_d$ in \cite{Zhao}, is that the subsets of ${\mathcal A}_d$ contain $0$. More precisely, $B_d=\{A\setminus\{0\}: A \in {\mathcal A}_d\}$. We chose not to include $0$ for simplifying the presentation.

Observe that,
if $A\in B_d$, then
\begin{equation}\label{boundcardinalityA}\#A\leq \left\lfloor\frac{d-1}{2}\right\rfloor
  \end{equation}
and
\begin{equation}\label{cardinalityB}\#B_d=3^{\lfloor\frac{d-1}{2}\rfloor}.\end{equation}
Indeed, in the case of $d$ even, $d/2$ can not belong to any $A\in B_d$. Furthermore, for both $d$ even and $d$ odd,
for each $A\in B_d$, for each $t\in[1,\lfloor\frac{d-1}{2}\rfloor]$, either (i) $t\in A$, $d-t\not\in A$, (ii) $t\not\in A$, $d-t\in A$, (iii) $t\not\in A$, $d-t\not\in A$.
To see \eqref{cardinalityB}, notice that, for any choice of (i) $t\in A$, $d-t\not\in A$, (ii) $t\not\in A$, $d-t\in A$, (iii) $t\not\in A$, $d-t\not\in A$, for any choice of $t\in[1,\lfloor\frac{d-1}{2}\rfloor]$, we obtain a set $A\in B_d$. 

\renewcommand\span{{\mathrm {span}}}
For $A\in B_d$ define
\begin{eqnarray*}
  \span_d(A)&=&(A\cup(A+A))\cap[1,d-1]\\
  \sigma_d(A)&=&\#\span(S)
\end{eqnarray*}

It holds
\begin{equation}\label{sigmabounds}0\leq \sigma_d(A)\leq d-1.\end{equation} The
lower bound is attained by $A=\emptyset$.
If $d$ is odd, the upper bound is attained by $A=[1,\frac{d-1}{2}]$.
If $d$ is even and $d\geq 8$, the upper bound is attained by $A=[1,\frac{d-4}{2}]\cup\{\frac{d+2}{2}\}$. Indeed, $[1,\dots,d-4]\in\span_d(A)$, since the whole interval can be obtained as sums of pairs in $[1,\frac{d-4}{2}]$. Furthermore, (i) $d-3$ is either $\frac{d+2}{2}$ (if $d=8$) or $\frac{d-8}{2}+\frac{d+2}{2}\in\span_d(A)$ (if $d\geq 10$); (ii) $d-2=\frac{d-6}{2}+\frac{d+2}{2}\in\span_d(A)$ (if $d\geq 8$); (iii) $d-1=\frac{d-4}{2}+\frac{d+2}{2}\in\span_d(A)$.

It will be used later that, if $d=F-2m$, and $A\in B_d$, whenever $d\geq 0$, 
\begin{equation}\label{2m+span}
  ((2m+A)\cup((m+A)+(m+A)))\cap[2m+1,F-1]=2m+\span(A).
  \end{equation}

Define now
\begin{eqnarray*}
  \gamma_{d,F}(A)&=&F-2-\#A-\sigma_d(A)\\
\nu_{d,m}(A)&=&m-2-\sigma_d(A)  
\end{eqnarray*}
From the definitions of $\gamma_{d,F}(A)$ and $\nu_{d,m}(A)$, and by \eqref{boundcardinalityA}, one can deduce that, if $F-2m=d$, for any $A\in B_d$,
whenever $d\geq 0$, \begin{equation}\label{nulegamma}
 \nu_{d,m}(A)\leq \gamma_{d,F}(A).\end{equation} 

\begin{table}[ht]
\begin{center}
  \begin{tabular}{|l|l|l|l|l|l|}
    \hline
$A$ & $\#A$  & $\span_1(A)$ & $\sigma_{1}(A)$ & $\gamma_{1,F}(A)$ & $\nu_{1,m}(A)$
      \\
      \hline
$\emptyset$ & 0 & $\emptyset$     &0 & $F-2$ & $m-2$ \\
\hline
  \end{tabular}\end{center}
    \caption{Sets in $B_1$}
    \label{t1}
\end{table}
~
\begin{table}[ht]
  \begin{center}
  \begin{tabular}{|l|l|l|l|l|l|}
    \hline
$A$ & $\#A$  & $\span_2(A)$ & $\sigma_{2}(A)$ & $\gamma_{2,F}(A)$ & $\nu_{2,m}(A)$
      \\
      \hline
$\emptyset$ & 0 & $\emptyset$     &0 & $F-2$ & $m-2$ \\
\hline
  \end{tabular}\end{center}
      \caption{Sets in $B_2$}
    \label{t2}
\end{table}
~
\begin{table}[ht]
\begin{center}
  \begin{tabular}{|l|l|l|l|l|l|}
    \hline
$A$ & $\#A$  & $\span_3(A)$ & $\sigma_{3}(A)$ & $\gamma_{3,F}(A)$ & $\nu_{3,m}(A)$
      \\
      \hline
      $\emptyset$ & 0 & $\emptyset$ &0 & $F-2$ & $m-2$ \\
      $\{1\}$     & 1 & $\{1,2\}$   &2 & $F-5$ & $m-4$ \\
      $\{2\}$     & 1 & $\{2\}$   &1 & $F-4$ & $m-3$ \\
      \hline
  \end{tabular}\end{center}
    \caption{Sets in $B_3$}
    \label{t3}
\end{table}
~
\begin{table}[ht]
\begin{center}
  \begin{tabular}{|l|l|l|l|l|l|}
    \hline
$A$ & $\#A$  & $\span_4(A)$ & $\sigma_{4}(A)$ & $\gamma_{4,F}(A)$ & $\nu_{4,m}(A)$
      \\
      \hline
      $\emptyset$ & 0 & $\emptyset$ &0 & $F-2$ & $m-2$ \\
      $\{1\}$     & 1 & $\{1,2\}$   &2 & $F-5$ & $m-4$ \\
      $\{3\}$     & 1 & $\{3\}$   &1 & $F-4$ & $m-3$ \\
      \hline
  \end{tabular}\end{center}
    \caption{Sets in $B_4$}
    \label{t4}
\end{table}
~
\begin{table}[ht]
  \begin{center}
  \begin{tabular}{|l|l|l|l|l|l|}
    \hline
$A$ & $\#A$  & $\span_5(A)$ & $\sigma_{5}(A)$ & $\gamma_{5,F}(A)$ & $\nu_{5,m}(A)$
      \\
      \hline
      $\emptyset$ & 0 & $\emptyset$     &0 & $F-2$ & $m-2$ \\
      $\{1\}$     & 1 & $\{1,2\}$       &2 & $F-5$ & $m-4$ \\
      $\{2\}$     & 1 & $\{2,4\}$       &2 & $F-5$ & $m-4$ \\
      $\{3\}$     & 1 & $\{3\}$         &1 & $F-4$ & $m-3$ \\
      $\{4\}$     & 1 & $\{4\}$         &1 & $F-4$ & $m-3$ \\
      $\{1,2\}$     & 2 & $\{1,2,3,4\}$ &4 & $F-8$ & $m-6$ \\
      $\{1,3\}$     & 2 & $\{1,2,3,4\}$ &4 & $F-8$ & $m-6$ \\
      $\{2,4\}$     & 2 & $\{2,4\}$     &2 & $F-6$ & $m-4$ \\
      $\{3,4\}$     & 2 & $\{3,4\}$     &2 & $F-6$ & $m-4$ \\
      \hline
  \end{tabular}\end{center}
    \caption{Sets in $B_5$}
    \label{t5}
\end{table}
~
\begin{table}[ht]
\begin{center}
  \begin{tabular}{|l|l|l|l|l|l|}
    \hline
$A$ & $\#A$  & $\span_6(A)$ & $\sigma_{6}(A)$ & $\gamma_{6,F}(A)$ & $\nu_{6,m}(A)$
      \\
      \hline
      $\emptyset$ & 0 & $\emptyset$     &0 & $F-2$ & $m-2$ \\
      $\{1\}$     & 1 & $\{1,2\}$       &2 & $F-5$ & $m-4$ \\
      $\{2\}$     & 1 & $\{2,4\}$       &2 & $F-5$ & $m-4$ \\
      $\{4\}$     & 1 & $\{4\}$         &1 & $F-4$ & $m-3$ \\
      $\{5\}$     & 1 & $\{5\}$         &1 & $F-4$ & $m-3$ \\
      $\{1,2\}$     & 2 & $\{1,2,3,4\}$ &4 & $F-8$ & $m-6$ \\
      $\{1,4\}$     & 2 & $\{1,2,4,5\}$ &4 & $F-8$ & $m-6$ \\
      $\{2,5\}$     & 2 & $\{2,4,5\}$   &3 & $F-7$ & $m-5$ \\
      $\{4,5\}$     & 2 & $\{4,5\}$     &2 & $F-6$ & $m-4$ \\
      \hline
  \end{tabular}\end{center}
    \caption{Sets in $B_6$}
    \label{t6}
\end{table}
~
\begin{table}[ht]
\begin{center}
  \begin{tabular}{|l|l|l|l|l|l|}
    \hline
$A$ & $\#A$  & $\span_7(A)$ & $\sigma_{7}(A)$ & $\gamma_{7,F}(A)$ & $\nu_{7,m}(A)$
      \\
      \hline
$\emptyset$ & 0 & $\emptyset$     &0 & $F-2$ & $m-2$ \\
$\{1\}$     & 1 & $\{1,2\}$       &2 & $F-5$ & $m-4$ \\
$\{2\}$     & 1 & $\{2,4\}$       &2 & $F-5$ & $m-4$ \\
$\{3\}$     & 1 & $\{3,6\}$       &2 & $F-5$ & $m-4$ \\
$\{4\}$     & 1 & $\{4\}$         &1 & $F-4$ & $m-3$ \\
$\{5\}$     & 1 & $\{5\}$         &1 & $F-4$ & $m-3$ \\
$\{6\}$     & 1 & $\{6\}$         &1 & $F-4$ & $m-3$ \\
$\{1,2\}$   & 2 & $\{1,2,3,4\}$   &4 & $F-8$ & $m-6$ \\
$\{1,3\}$   & 2 & $\{1,2,3,4,6\}$ &5 & $F-9$ & $m-7$ \\
$\{1,4\}$   & 2 & $\{1,2,4,5\}$   &4 & $F-8$ & $m-6$ \\
$\{1,5\}$   & 2 & $\{1,2,5,6\}$   &4 & $F-8$ & $m-6$ \\
$\{2,3\}$   & 2 & $\{2,3,4,5,6\}$ &5 & $F-9$ & $m-7$ \\
$\{2,4\}$   & 2 & $\{2,4,6\}$     &3 & $F-7$ & $m-5$ \\
$\{2,6\}$   & 2 & $\{2,4,6\}$     &3 & $F-7$ & $m-5$ \\
$\{3,5\}$   & 2 & $\{3,5,6\}$     &3 & $F-7$ & $m-5$ \\
$\{3,6\}$   & 2 & $\{3,6\}$       &2 & $F-6$ & $m-4$ \\
$\{4,5\}$   & 2 & $\{4,5\}$       &2 & $F-6$ & $m-4$ \\
$\{4,6\}$   & 2 & $\{4,6\}$       &2 & $F-6$ & $m-4$ \\
$\{5,6\}$   & 2 & $\{5,6\}$       &2 & $F-6$ & $m-4$ \\
$\{1,2,3\}$ & 3 & $\{1,2,3,4,5,6\}$ &6 & $F-11$ & $m-8$ \\
$\{1,2,4\}$ & 3 & $\{1,2,3,4,5,6\}$ &6 & $F-11$ & $m-8$ \\
$\{1,3,5\}$ & 3 & $\{1,2,3,4,5,6\}$ &6 & $F-11$ & $m-8$ \\
$\{1,4,5\}$ & 3 & $\{1,2,4,5,6\}$   &5 & $F-10$ & $m-7$ \\
$\{2,3,6\}$ & 3 & $\{2,3,4,5,6\}$   &5 & $F-10$ & $m-7$ \\
$\{2,4,6\}$ & 3 & $\{2,4,6\}$       &3 & $F-8$ & $m-5$ \\
$\{3,5,6\}$ & 3 & $\{3,5,6\}$       &3 & $F-8$ & $m-5$ \\
$\{4,5,6\}$ & 3 & $\{4,5,6\}$       &3 & $F-8$ & $m-5$ \\
\hline
  \end{tabular}\end{center}
    \caption{Sets in $B_7$}
    \label{t7}
\end{table}
~
\begin{table}[ht]
\begin{center}
  \begin{tabular}{|l|l|l|l|l|l|}
    \hline
$A$ & $\#A$  & $\span_8(A)$ & $\sigma_{8}(A)$ & $\gamma_{8,F}(A)$ & $\nu_{8,m}(A)$
      \\
      \hline
$\emptyset$ & 0&$\emptyset$      & 0& $F-2$ & $m-2$ \\
$\{1\}$     & 1& $\{1,2\}$       & 2& $F-5$ & $m-4$ \\
$\{2\}$     & 1& $\{2,4\}$       & 2& $F-5$ & $m-4$ \\
$\{3\}$     & 1& $\{3,6\}$       & 2& $F-5$ & $m-4$ \\
$\{5\}$     & 1& $\{5\}$         & 1& $F-4$ & $m-3$ \\
$\{6\}$     & 1& $\{6\}$         & 1& $F-4$ & $m-3$ \\
$\{7\}$     & 1& $\{7\}$         & 1& $F-4$ & $m-3$ \\
$\{1,2\}$   & 2& $\{1,2,3,4\}$   & 4& $F-8$ & $m-6$ \\
$\{1,3\}$   & 2& $\{1,2,3,4,6\}$ & 5& $F-9$ & $m-7$ \\
$\{1,5\}$   & 2& $\{1,2,5,6\}$   & 4& $F-8$ & $m-6$ \\
$\{1,6\}$   & 2& $\{1,2,6,7\}$   & 4& $F-8$ & $m-6$ \\
$\{2,3\}$   & 2& $\{2,3,4,5,6\}$ & 5& $F-9$ & $m-7$ \\
$\{2,5\}$   & 2& $\{2,4,5,7\}$   & 4& $F-8$ & $m-6$ \\
$\{2,7\}$   & 2& $\{2,4,7\}$     & 3& $F-7$ & $m-5$ \\
$\{3,6\}$   & 2& $\{3,6\}$       & 2& $F-6$ & $m-4$ \\
$\{3,7\}$   & 2& $\{3,6,7\}$     & 3& $F-7$ & $m-5$ \\
$\{5,6\}$   & 2& $\{5,6\}$       & 2& $F-6$ & $m-4$ \\
$\{5,7\}$   & 2& $\{5,7\}$       & 2& $F-6$ & $m-4$ \\
$\{6,7\}$   & 2& $\{6,7\}$       & 2& $F-6$ & $m-4$ \\
$\{1,2,3\}$ & 3& $\{1,2,3,4,5,6\}$   & 6& $F-11$ & $m-8$ \\
$\{1,2,5\}$ & 3& $\{1,2,3,4,5,6,7\}$ & 7& $F-12$ & $m-9$ \\
$\{1,3,6\}$ & 3& $\{1,2,3,4,6,7\}$   & 6& $F-11$ & $m-8$ \\
$\{1,5,6\}$ & 3& $\{1,2,5,6,7\}$     & 5& $F-10$ & $m-7$ \\
$\{2,3,7\}$ & 3& $\{2,3,4,5,6,7\}$   & 6& $F-11$ & $m-8$ \\
$\{2,5,7\}$ & 3& $\{2,4,5,7\}$       & 4& $F-9$ & $m-6$ \\
$\{3,6,7\}$ & 3& $\{3,6,7\}$         & 3& $F-8$ & $m-5$ \\
      $\{5,6,7\}$ & 3& $\{5,6,7\}$         & 3& $F-8$ & $m-5$ \\
      \hline
    \end{tabular}\end{center}
    \caption{Sets in $B_8$}
    \label{t8}
\end{table}

To illustrate all these definitions, in Table~\ref{t1} to Table~\ref{t8} we provide the list of all sets in $B_1,B_2,\dots,B_8$, respectively, with their associated parameters. As explained, in all cases there are $3^{\lfloor\frac{d-1}{2}\rfloor}$ sets, that is, $\#B_1=\#B_2=1$, $\#B_3=\#B_4=3$, $\#B_5=\#B_6=9$.

\begin{lemma}\label{nubound}
  Let $F,m$ be nonnegative integers such that $\frac{F+1}{3}\leq m\leq \frac{F}{2}$,
  and let $d=F-2m$. Then, for any $A\in B_d$, it holds 
$\nu_{d,m}(A)\geq 0$.
\end{lemma}

\begin{proof}
  We have $\nu_{d,m}(A)=m-2-\sigma_d(A)$ which, by \eqref{sigmabounds}, is at least
  $3m-F-1$, and which, by the hypothesis, is at least $0$.
  \end{proof}

For any given integers $F,m$ ($F\geq 2m$) and for any $A\in B_{F-2m}$, define
$$M_{F,m}(A)=\{0,m\}\cup(m+A)\cup\{2m\}\cup(2m+\span_{F-2m}(A))\cup[F+1,\infty).$$

As an example, take $F=28$, $m=10$. Then $F-2m=8$ as in Table~\ref{t8}. Consider $A=\{2,5\}$. Then, as shown in Table~\ref{t8}, $\span_8(\{2,5\})=\{2,4,5,7\}$, and
  $$M_{28,10}(\{2,5\})=\{0,10\}\cup \{12,15\} \cup \{20\} \cup \{22,24,25,27\}\cup [29,\infty)$$

\begin{lemma}\label{MAVT}
  Let $F,m$ be nonnegative integers such that $\frac{F+1}{3}\leq m\leq \frac{F}{2}$,
  and let $d=F-2m$. For any $A\in B_d$, the following statements hold.
  \begin{enumerate}
  \item The sets in the union in the definition of $M_{F,m}(A)$ are disjoint.
  \item $M_{F,m}(A)$ is a numerical semigroup.
  \item The genus of $M_{F,m}(A)$ is $\gamma_{d,F}(A)$.
  \item Any numerical semigroup in ${\mathcal S}_{F,m,g}$ containing $m+A$, contains $M_{F,m}(A)$.
  \end{enumerate}
\end{lemma}

\begin{proof}
  \begin{enumerate}
  \item It is a consequence of the inequalities
    \begin{itemize}
    \item $m<\min(m+A)$,
    \item $\max(m+A)\leq m+d-1=F-m-1\leq 2m-2$ (where we used the hypothesis of the lemma),
    \item $2m<\min(2m+\span_d(A))$,
      \item $\max(2m+\span_d(A))\leq 2m+d-1=F-1$.
    \end{itemize}
  \item It is enough to prove that it is closed under addition.
    Indeed, one can check that
    \begin{itemize}
    \item $m+M_{F,m}(A)\subseteq M_{F,m}(A)$
    \item $(m+A)+(m+A)\subseteq 2m+(A+A)=\left((2m+(A+A))\cap[1,F-1]\right)\cup\left((2m+(A+A))\cap[F,\infty)\right)$,
      But,
      \begin{itemize}
      \item $(2m+(A+A))\cap[1,F-1]\subseteq 2m+\span_d(A)\subseteq M_{F,m}(A)$, by \eqref{2m+span},
      \item  $(2m+(A+A))\cap[F,\infty)\subset[F+1,\infty)\subseteq M_{F,m}(A)$, because $F\not\in 2m+(A+A)$, by definition of $B_d$.
      \end{itemize}
        \item All elements in $(m+A)+(\{2m\}\cup(2m+\span_d(S)))$ are larger than or equal to $3m+1\geq F+1$, and, so, they belong to $M_{F,m}(A)$.
          \item All the sums of two elements belonging to $\{2m\}\cup(2m+\span_d(A))\cup [F+1,\infty)$ will be larger than or equal to $4m\geq F+1$, and so, all these sums belong to $M_{F,m}(A)$.
    \end{itemize}
  \item By definition of $M_{F,m}(A)$, its genus is
    \begin{eqnarray*}g(M_{F,m}(A))&=&F-\#\left(\{m\}\cup(m+A)\cup\{2m\}\cup(2m+\span_d(A))\right)\\&=&F-1-\#A-1-\sigma_d(A)=\gamma_{d,F}(A).\end{eqnarray*}
    \item Any numerical semigroup in ${\mathcal S}_{F,m,g}$ containing $m+A$ will contain $\{0,m\}$, $m+A$, $\{2m\}$, $[F+1,\infty)$, and the union of $m+(m+A)=2m+A$ and $(m+A)+(m+A)$. If it contains the union of $2m+A$ and $(m+A)+(m+A)$, by \eqref{2m+span}, then it needs to contain $2m+\span_d(A)$. Hence, any numerical semigroup in ${\mathcal S}_{F,m,g}$ containing $m+A$, contains $M_{F,m}(A)$.
    \end{enumerate}
\end{proof}

One can check all these statements in the case of $M_{28,10}(\{2,5\})$. In particular, its genus is $20$, which coincides with $\gamma_{d,F}(\{2,5\})$, which, as predicted in Table~\ref{t8}, is $F-8$.

Now, consider a general semigroup $S$ of Frobenius number $F$, multiplicity $m$ and genus equal to $g$, that is, $S\in {\mathcal S}_{F,m,g}$. For instance take in ${\mathcal S}_{28,10,18}$ the semigroup
$$S=\{0,10,12,15,19,20,22,23,24,25,27\}\cup[29,\infty).$$
The next lemma proves that $S\cap[m+1,F-m-1]=m+A$ for some $A$ in $B_d$. In our example, $S\cap[11,17]=\{12,15\}=10+\{2,5\}$, with $\{2,5\}\in B_4.$

\begin{lemma}\label{AdeS}
  Let $F,m,g$ be nonnegative integers such that $\frac{F+1}{3}\leq m\leq \frac{F}{2}$, and let $d=F-2m$.   
  Let $S\in{\mathcal S}_{F,m,g}$, and let $A$ be such that 
  $S\cap[m+1,F-m-1]=m+A$. That is, let $A=-m+(S\cap[m+1,F-m-1])$. Then,
  \begin{enumerate}
  \item $A\in B_d$,
  \item $\gamma_{d,F}(A)\geq g$,
  \item $\nu_{d,m}(A)\geq \gamma_{d,F}(A)-g$.
    \end{enumerate}
\end{lemma}

\begin{proof}
  \begin{enumerate}
    \item
  On one hand, $-m+(S\cap[m+1,F-m-1])\subseteq [1,F-2m-1]$. On the other hand, if $t$ and $F-2m-t$ both belong to $-m+(S\cap[m+1,F-m-1])$, then $t+m$ and $F-m-t$ both belong to $S\cap[m+1,F-m-1]$. Hence, $(t+m)+(F-m-t)=F$ belongs to $S$, which is a contradiction because $F$ is the Frobenius number of $S$.
\item 
  The inequality $\gamma_{d,F}(A)\geq g$ follows from the fact that $\gamma_{d,F}(A)$ and $g$ are the genera of $M_{F,m}(A)$ and $S$, respectively, while $M_{F,m}(A)\subseteq S$, by Lemma~\ref{MAVT}.
\item Observe that $\nu_{d,m}(A)=\#\left([F-m+1,F-1]\setminus M_{F,m}(A)\right)$. Indeed,
  $$\#\left([F-m+1,F-1]\cap M_{F,m}(A)\right)=\#\{2m\}+\#\span_d(A)=1+\sigma_d(A)$$
  and, so,
  $$\#\left([F-m+1,F-1]\setminus M_{F,m}(A)\right)=m-1-1-\sigma_d(A)=\nu_{d,m}(A).$$
  Now,
  $$S\setminus M_{F,m}(A)\subseteq [2m,F-1]\setminus M_{F,m}(A) \subseteq [F-m+1,F-1]\setminus M_{F,m}(A)$$
  and, so,
by Lemma~\ref{MAVT},
  $$\gamma_{d,F}(A)-g=\#\left(S\setminus M_{F,m}(A)\right)\leq \#\left([F-m+1,F-1]\setminus M_{F,m}(A)\right)=\nu_{d,m}(A).$$
  
  \end{enumerate}
  \end{proof}

The previous lemma suggests defining, for nonnegative integers $F,m,g$ with $F-2m\geq 0$, and for any $A\in B_{F-2m}$, $${\mathcal S}_{F,m,g}^A=\{S\in{\mathcal S}_{F,m,g}\mbox{ such that }S\cap[m+1,F-m-1]=m+A\}.$$
Now, the next corollary is a consequence of Lemma~\ref{AdeS}.

\begin{corollary}\label{cor}
  If $\frac{F+1}{3}\leq m\leq \frac{F}{2}$, 
  then
  \begin{eqnarray}\label{eqSfmg}{\mathcal S}_{F,m,g}&=&\bigsqcup_{A\in B_{d}}{\mathcal S}_{F,m,g}^A,\\
          {\mathcal S}_{F,m}&=&\bigsqcup_{A\in B_{d}}\left(\bigsqcup_{g=\gamma_{d,F}(A)-\nu_{d,m}(A)}^{\gamma_{d,F}(A)}{\mathcal S}_{F,m,g}^A\right).\label{eqSfm}\\
    \end{eqnarray}
  where $d=F-2m$.
\end{corollary}

\section{Counting numerical semigroups with $m\geq\frac{F+1}{3}$}
\label{s:counting}

Define
$${\mathcal G}_{F,m,g}^A=\{G \mbox{ subset of }{[F-m+1,F-1]\setminus(\{2m\}\cup(2m+\span_d(A)))}\mbox{ such that }\#G=\gamma_{d,F}(A)-g\},$$
where, as before, $d=F-2m$.
\begin{lemma}\label{SG}
    Let $F,m,g$ be nonnegative integers such that $\frac{F+1}{3}\leq m\leq \frac{F}{2}$, and let $d=F-2m$. Let $A\in B_d$. Then,
  \begin{enumerate}
  \item ${\mathcal S}_{F,m,g}^A$ is in bijection with ${\mathcal G}_{F,m,g}^A$,
  \item $\#{\mathcal G}_{F,m,g}^A=\binom{\nu_{d,m}(A)}{\gamma_{d,F}(A)-g}$, which is well defined by Corollary~\ref{cor},
  \end{enumerate}
and, therefore, $\#{\mathcal S}_{F,m,g}^A=\binom{\nu_{d,m}(A)}{\gamma_{d,F}(A)-g}$. 
\end{lemma}

\begin{proof}
  \begin{enumerate}
  \item
Define the map $$\begin{array}{rcl}\mu:{\mathcal S}_{F,m,g}^A&\longrightarrow&{\mathcal G}_{F,m,g}^A\\S&\mapsto&S\setminus M_{F,m}(A)\end{array}.$$
The map $\mu$ is injective. Indeed, if $S,S'\in{\mathcal S}_{F,m,g}^A$ and $S\neq S'$, then, since $M_{F,m}(A)\subseteq S,S'$, then $S\setminus M_{F,m}(A)\neq S'\setminus M_{F,m}(A)$.

Let us see that $\mu$ is exhaustive. Let us first see that, if $G\in{\mathcal G}_{F,m,g}^A$, then $M_{F,m}(A)\cup G$ is a numerical semigroup. Indeed,
on one hand, $G+G\subseteq M_{F,m}(A)\cup G$ because $\min(G+G)=2\min G\geq 2F-2m+2\geq F+1$, by the hypothesis on $m$.
    Let us see, on the other hand, that $M_{F,m}(A)+G\subseteq M_{F,m}(A)\cup G$. For $a\in M_{F,m}(A)$, it can be either $a=0$ or $a\geq m$. If $a=0$, then $a+G=G\subseteq M_{F,m}(A)\cup G$. Otherwise, $a+G\subseteq [a+\min G,\infty)\subseteq [F+1,\infty)\subseteq M_{F,m}(A)\cup G$. This concludes that $M_{F,m}(A)\cup G\in{\mathcal S}_{F,m,g}^A$. Exhaustivity follows now from the fact that $M_{F,m}(A)\cap G=\emptyset$ and, so, $\mu(M_{F,m}(A)\cup G)=G$. 
  
  \item Since $\{2m\}\cup(2m+\span_d(A))\subseteq [F-m+1,F-1]$, the cardinality of $[F-m+1,F-1]\setminus(\{2m\}\cup(2m+\span_d(A)))$ is $m-2-\sigma_d(A)=\nu_{d,m}(A)$.
    The result follows immediately. 
    \end{enumerate}
  \end{proof}

In the next theorem, we assume that $\binom{a}{b}=0$ whenever $b<0$ or $b>a$.

\begin{theorem}\label{pot2}
  If $\frac{F+1}{3}\leq m\leq \frac{F}{2}$, and, as ususal, $\frac{F+1}{2}\leq g\leq F$,
  then
  $$n_{F,m,g}=\sum_{A\in B_{d}}\binom{\nu_{d,m}(A)}{\gamma_{d,F}(A)-g}$$
and  
  $$n_{F,m}=\sum_{A\in B_d}2^{\nu_{d,m}(A)},$$
  where $d=F-2m$.
  \end{theorem}

\begin{proof}
  The first formula is a direct consequence of
\eqref{eqSfmg} in
  Corollary~\ref{cor} and Lemma~\ref{SG}.
  As for the second formula, 
by \eqref{eqSfm},
  $$n_{F,m}=\sum_{A\in B_{d}}\sum_{g=\gamma_{d,F}(A)-\nu_{d,m}(A)}^{\gamma_{d,F}(A)}\binom{\nu_{d,m}(A)}{\gamma_{d,F}(A)-g}=\sum_{A\in B_{d}}\sum_{k=0}^{\nu_{d,m}(A)}\binom{\nu_{d,m}(A)}{k}=\sum_{A\in B_d}2^{\nu_{d,m}(A)}.$$
\end{proof}

As an example of the first formula, take $F=24$, $m=10$. In this case, $d=4$, and, considering the list of sets in $B_4$ shown in Table~\ref{t4}, 
$$n_{24,10,g}=\binom{m-2}{F-2-g}+\binom{m-4}{F-5-g}+\binom{m-3}{F-4-g}=\binom{8}{22-g}+\binom{6}{19-g}+\binom{7}{20-g}$$
Now, evaluating for instance at $g=15$ we obtain
$$n_{24,10,15}=\binom{8}{7}+\binom{6}{4}+\binom{7}{5}=8+15+21=44,$$
as can be checked in Table~\ref{tF24}.

As an example of the second formula, considering $F=24$ and $m=10$, and looking at Table~\ref{t4}, we get that, $n_{F,m}=2^{m-2}+2^{m-4}+2^{m-3}=2^8+2^6+2^7=256+64+128=448$, as can be cheched in Table~\ref{tF24}.

As a larger example, looking at Table~\ref{t7}, we get that, for any $m\geq 8$, and for $F=2m+7$,
\begin{eqnarray*}
  n_{F,m}&=&1\cdot 2^{m-2}+    3\cdot 2^{m-3}+    7\cdot 2^{m-4}+    6\cdot 2^{m-5}+    3\cdot 2^{m-6}+    4\cdot 2^{m-7}+    3\cdot 2^{m-8} \\&=& 2^{m-8} \cdot (64 + 3 \cdot 32 + 7\cdot 16+6\cdot 8+3\cdot 4+4\cdot 2+3\cdot 1) \\&=& 343 \cdot 2^{m-8}.
\end{eqnarray*}

On the other hand, looking at Table~\ref{t8}, we get that, for any $m\geq 9$, for $F=2m+8$,
\begin{eqnarray*}
 n_{F,m}&=&1\cdot 2^{m-2}+    3\cdot 2^{m-3}+    7\cdot 2^{m-4}+    4\cdot 2^{m-5}+   5\cdot 2^{m-6}+    3\cdot 2^{m-7}+    3\cdot 2^{m-8} + 1\cdot 2^{m-9}
  \\&=& 2^{m-9} \cdot (128 + 3 \cdot 64 + 7 \cdot 32 + 4\cdot 16+5\cdot 8+3\cdot 4+3\cdot 2+ 1) \\&=& 667 \cdot 2^{m-9}.
\end{eqnarray*}

Similarly, the same argument proves the formulas
$$\begin{array}{rcll}
n_{2m,m}&=&0&\\
n_{2m+1,m}&=&2^{m-2}&\mbox{if }m\geq 2\\
n_{2m+2,m}&=&2^{m-2}&\mbox{if }m\geq 3\\
n_{2m+3,m}&=&7\cdot 2^{m-4}&\mbox{if }m\geq 4\\
n_{2m+4,m}&=&7\cdot 2^{m-4}&\mbox{if }m\geq 5\\
n_{2m+5,m}&=&50\cdot 2^{m-6}&\mbox{if }m\geq 6\\
n_{2m+6,m}&=&96\cdot 2^{m-7}&\mbox{if }m\geq 7\\
n_{2m+7,m}&=&343\cdot 2^{m-8}&\mbox{if }m\geq 8\\
n_{2m+8,m}&=&667\cdot 2^{m-9}&\mbox{if }m\geq 9\\
&\vdots&&\\
\end{array}$$

These formulas can be checked in Table~\ref{tF25} and Table~\ref{tF26}, since one can check that
$$\begin{array}{rcll}
n_{26,13}&=&0&\\
n_{25,12}&=&2^{m-2}=2^{10}=1024&\\
n_{26,12}&=&2^{m-2}=2^{10}=1024&\\
n_{25,11}&=&7\cdot 2^{m-4}=7\cdot 2^7=896 &\\
n_{26,11}&=&7\cdot 2^{m-4}=7\cdot 2^7=896&\\
n_{25,10}&=&50\cdot 2^{m-6}=50\cdot 2^4=800&\\
n_{26,10}&=&96\cdot 2^{m-7}=96\cdot 2^3=768&\\
n_{25,9}&=&343\cdot 2^{m-8}=343\cdot 2=686&\\
n_{26,9}&=&667\cdot 2^{m-9}=667&\\
\end{array}$$

\section{Two-step-doubling growth with respect to the Frobenius number}
\label{s:Fgrowth}

Now we want to relate $n_{F,m}$ to $n_{F-2,m-1}.$
The next lemma gives a well known formula for $n_{F,m}$ for the case in which $m>\frac{F}{2}.$

\begin{lemma}\label{mgran}
    Let $F,m$ be nonnegative integers  with $m> \frac{F}{2}$. Then, $$n_{F,m}=2^{F-m-1}.$$
  \end{lemma}

\begin{proof}
If $m> \frac{F}{2}$, the semigroups with Frobenius number $F$ and multiplicity $m$ are in bijection with the subsets of $[m+1,F-1]$, where the correspondence is given by assigning to a semigroup the subset of nongaps between $m+1$ and $F-1$. This proves the formula.
  \end{proof}

\begin{theorem}\label{t:twostepdoubling}
  Let $F,m$ be nonnegative integers. If $\frac{F+2}{3}\leq m<F-1$, then $$n_{F,m}=2n_{F-2,m-1}.$$
\end{theorem}

\begin{proof}
  Suppose first that $m>\frac{F}{2}$. By Lemma~\ref{mgran}, $$n_{F,m}=2^{F-m-1}.$$
  But, if $m>\frac{F}{2}$, we also have $m-1>\frac{F-2}{2}$ and, again, by Lemma~\ref{mgran} applied to $F-2$ and $m-1$, we get $$n_{F-2,m-1}=2^{F-m-2}.$$ Then, it follows that
  $$n_{F,m}=2n_{F-2,m-1}.$$
  Suppose now that $\frac{F+2}{3}\leq m\leq \frac{F}{2}$. By Theorem~\ref{pot2},
  $$n_{F,m}=\sum_{A\in B_d}2^{\nu_{d,m}(A)},$$
  where $d=F-2m$.
  But, if $\frac{F+2}{3}\leq m\leq \frac{F}{2}$, we also have $\frac{(F-2)+1}{3}\leq m-1\leq \frac{F-2}{2}$ and, again, by Theorem~\ref{pot2} applied to $F-2$ and $m-1$, we get $$n_{F-2,m-1}=\sum_{A\in B_d}2^{\nu_{d,m-1}(A)},$$
  where $d=(F-2)-2(m-1)=F-2m$ stays the same.

   Notice that $\nu_{d,m}(A)=m-2-\sigma_d(A)=1+(m-1)-2-\sigma_d(A)=1+\nu_{d,m-1}(A)$.
   Hence,
$$n_{F,m}=\sum_{A\in B_d}2^{\nu_{d,m}(A)}=\sum_{A\in B_d}2^{1+\nu_{d,m-1}(A)}=2\sum_{A\in B_d}2^{\nu_{d,m-1}(A)}=2 n_{F-2,m-1}.$$
  \end{proof}

The previous theorem can be checked looking at the last columns of
Table~\ref{tF24} and Table~\ref{tF26}.
For instance, $n_{26,10}=768$, which coincides with $2\cdot n_{24,9}=2\cdot 384$.
Table~\ref{taulaFm} gives all values of $n_{F,m}$ for $F$ up to $30$. The values where Theorem~\ref{t:twostepdoubling} holds are highlighted in bold.

We want to remark that the hypotheses of the theorem are necessary. Indeed,
(i) if $m=F-1$, then $n_{F,F-1}=1$, while $n_{F-2,F-2}=0$;
(ii) if $m=F$, then $n_{F,F}=0$, while $n_{F-2,F-1}=1$;
(iii) if $m=F+1$, then $n_{F,F+1}=1$, while $n_{F-2,F}=0$.
On the other hand, we find cases with $m<\frac{F+2}{3}$ for which the equality is not satisfied,
for instance:
(i) $n_{21,5}=31$, while $n_{19,4}=16$; (ii) $n_{21,4}=17$; while $n_{19,3}=4$. This is just to show different examples in which $n_{F,m}>2n_{F-2,m-1}$ and in which $n_{F,m}<2n_{F-2,m-1}$.

  \begin{table}
  \rotatebox{90}{\resizebox{.8\textheight}{!}{{\renewcommand{\arraystretch}{2}\begin{tabular}{||l|rrrrrrrrrrrrrrrrrrrrrrrrrrrrrr|r||}\hline\hline $F\setminus m$   & m=2  & m=3  & m=4  & m=5  & m=6  & m=7  & m=8  & m=9  & m=10  & m=11  & m=12  & m=13  & m=14  & m=15  & m=16  & m=17  & m=18  & m=19  & m=20  & m=21  & m=22  & m=23  & m=24  & m=25  & m=26  & m=27  & m=28  & m=29  & m=30  & m=31 & total \\\hline\hline
F=1 & 1 &&&&&&&&&&&&&&&&&&&&&&&&&&&&&& 1 \\\hline
 F=2 & 0 & 1 &&&&&&&&&&&&&&&&&&&&&&&&&&&&& 1 \\\hline
 F=3 & 1 & 0 & 1 &&&&&&&&&&&&&&&&&&&&&&&&&&&& 2 \\\hline
 F=4 & 0 & 1 & 0 & 1 &&&&&&&&&&&&&&&&&&&&&&&&&&& 2 \\\hline
 F=5 & 1 & 2 & 1 & 0 & 1 &&&&&&&&&&&&&&&&&&&&&&&&&& 5 \\\hline
 F=6 & 0 & 0 & 2 & 1 & 0 & 1 &&&&&&&&&&&&&&&&&&&&&&&&& 4 \\\hline
 F=7  & 1  & {\bf 2}  & {\bf 4}  & {\bf 2}  & 1  & 0  & 1 &&&&&&&&&&&&&&&&&&&&&&&& 11 \\\hline
F=8  & 0  & {\bf 2}  & {\bf 0}  & {\bf 4}  & {\bf 2}  & 1  & 0  & 1 &&&&&&&&&&&&&&&&&&&&&&& 10 \\\hline
F=9  & 1  & {\bf 0}  & {\bf 4}  & {\bf 8}  & {\bf 4}  & {\bf 2}  & 1  & 0  & 1 &&&&&&&&&&&&&&&&&&&&&& 21 \\\hline
F=10  & 0  & 2  & {\bf 4}  & {\bf 0}  & {\bf 8}  & {\bf 4}  & {\bf 2}  & 1  & 0  & 1 &&&&&&&&&&&&&&&&&&&&& 22 \\\hline
F=11  & 1  & 3  & {\bf 7}  & {\bf 8}  & {\bf 16}  & {\bf 8}  & {\bf 4}  & {\bf 2}  & 1  & 0  & 1 &&&&&&&&&&&&&&&&&&&& 51 \\\hline
F=12  & 0  & 0  & {\bf 0}  & {\bf 8}  & {\bf 0}  & {\bf 16}  & {\bf 8}  & {\bf 4}  & {\bf 2}  & 1  & 0  & 1 &&&&&&&&&&&&&&&&&&& 40 \\\hline
F=13  & 1  & 3  & 8  & {\bf 14}  & {\bf 16}  & {\bf 32}  & {\bf 16}  & {\bf 8}  & {\bf 4}  & {\bf 2}  & 1  & 0  & 1 &&&&&&&&&&&&&&&&&& 106 \\\hline
F=14  & 0  & 3  & 6  & {\bf 14}  & {\bf 16}  & {\bf 0}  & {\bf 32}  & {\bf 16}  & {\bf 8}  & {\bf 4}  & {\bf 2}  & 1  & 0  & 1 &&&&&&&&&&&&&&&&& 103 \\\hline
F=15  & 1  & 0  & 11  & {\bf 0}  & {\bf 28}  & {\bf 32}  & {\bf 64}  & {\bf 32}  & {\bf 16}  & {\bf 8}  & {\bf 4}  & {\bf 2}  & 1  & 0  & 1 &&&&&&&&&&&&&&&& 200 \\\hline
F=16  & 0  & 3  & 0  & 14  & {\bf 28}  & {\bf 32}  & {\bf 0}  & {\bf 64}  & {\bf 32}  & {\bf 16}  & {\bf 8}  & {\bf 4}  & {\bf 2}  & 1  & 0  & 1 &&&&&&&&&&&&&&& 205 \\\hline
F=17  & 1  & 4  & 12  & 22  & {\bf 50}  & {\bf 56}  & {\bf 64}  & {\bf 128}  & {\bf 64}  & {\bf 32}  & {\bf 16}  & {\bf 8}  & {\bf 4}  & {\bf 2}  & 1  & 0  & 1 &&&&&&&&&&&&&& 465 \\\hline
F=18  & 0  & 0  & 9  & 20  & {\bf 0}  & {\bf 56}  & {\bf 64}  & {\bf 0}  & {\bf 128}  & {\bf 64}  & {\bf 32}  & {\bf 16}  & {\bf 8}  & {\bf 4}  & {\bf 2}  & 1  & 0  & 1 &&&&&&&&&&&&& 405 \\\hline
F=19  & 1  & 4  & 16  & 31  & 57  & {\bf 100}  & {\bf 112}  & {\bf 128}  & {\bf 256}  & {\bf 128}  & {\bf 64}  & {\bf 32}  & {\bf 16}  & {\bf 8}  & {\bf 4}  & {\bf 2}  & 1  & 0  & 1 &&&&&&&&&&&& 961 \\\hline
F=20  & 0  & 4  & 0  & 0  & 48  & {\bf 96}  & {\bf 112}  & {\bf 128}  & {\bf 0}  & {\bf 256}  & {\bf 128}  & {\bf 64}  & {\bf 32}  & {\bf 16}  & {\bf 8}  & {\bf 4}  & {\bf 2}  & 1  & 0  & 1 &&&&&&&&&&& 900 \\\hline
F=21  & 1  & 0  & 17  & 31  & 75  & {\bf 0}  & {\bf 200}  & {\bf 224}  & {\bf 256}  & {\bf 512}  & {\bf 256}  & {\bf 128}  & {\bf 64}  & {\bf 32}  & {\bf 16}  & {\bf 8}  & {\bf 4}  & {\bf 2}  & 1  & 0  & 1 &&&&&&&&&& 1828 \\\hline
F=22  & 0  & 4  & 12  & 32  & 68  & 101  & {\bf 192}  & {\bf 224}  & {\bf 256}  & {\bf 0}  & {\bf 512}  & {\bf 256}  & {\bf 128}  & {\bf 64}  & {\bf 32}  & {\bf 16}  & {\bf 8}  & {\bf 4}  & {\bf 2}  & 1  & 0  & 1 &&&&&&&&& 1913 \\\hline
F=23  & 1  & 5  & 22  & 43  & 122  & 152  & {\bf 343}  & {\bf 400}  & {\bf 448}  & {\bf 512}  & {\bf 1024}  & {\bf 512}  & {\bf 256}  & {\bf 128}  & {\bf 64}  & {\bf 32}  & {\bf 16}  & {\bf 8}  & {\bf 4}  & {\bf 2}  & 1  & 0  & 1 &&&&&&&& 4096 \\\hline
F=24  & 0  & 0  & 0  & 42  & 0  & 144  & {\bf 0}  & {\bf 384}  & {\bf 448}  & {\bf 512}  & {\bf 0}  & {\bf 1024}  & {\bf 512}  & {\bf 256}  & {\bf 128}  & {\bf 64}  & {\bf 32}  & {\bf 16}  & {\bf 8}  & {\bf 4}  & {\bf 2}  & 1  & 0  & 1 &&&&&&& 3578 \\\hline
F=25  & 1  & 5  & 23  & 0  & 132  & 228  & 382  & {\bf 686}  & {\bf 800}  & {\bf 896}  & {\bf 1024}  & {\bf 2048}  & {\bf 1024}  & {\bf 512}  & {\bf 256}  & {\bf 128}  & {\bf 64}  & {\bf 32}  & {\bf 16}  & {\bf 8}  & {\bf 4}  & {\bf 2}  & 1  & 0  & 1 &&&&&& 8273 \\\hline
F=26  & 0  & 5  & 16  & 44  & 109  & 216  & 334  & {\bf 667}  & {\bf 768}  & {\bf 896}  & {\bf 1024}  & {\bf 0}  & {\bf 2048}  & {\bf 1024}  & {\bf 512}  & {\bf 256}  & {\bf 128}  & {\bf 64}  & {\bf 32}  & {\bf 16}  & {\bf 8}  & {\bf 4}  & {\bf 2}  & 1  & 0  & 1 &&&&& 8175 \\\hline
F=27  & 1  & 0  & 28  & 58  & 161  & 319  & 561  & {\bf 0}  & {\bf 1372}  & {\bf 1600}  & {\bf 1792}  & {\bf 2048}  & {\bf 4096}  & {\bf 2048}  & {\bf 1024}  & {\bf 512}  & {\bf 256}  & {\bf 128}  & {\bf 64}  & {\bf 32}  & {\bf 16}  & {\bf 8}  & {\bf 4}  & {\bf 2}  & 1  & 0  & 1 &&&& 16132 \\\hline
F=28  & 0  & 5  & 0  & 57  & 142  & 0  & 450  & 711  & {\bf 1334}  & {\bf 1536}  & {\bf 1792}  & {\bf 2048}  & {\bf 0}  & {\bf 4096}  & {\bf 2048}  & {\bf 1024}  & {\bf 512}  & {\bf 256}  & {\bf 128}  & {\bf 64}  & {\bf 32}  & {\bf 16}  & {\bf 8}  & {\bf 4}  & {\bf 2}  & 1  & 0  & 1 &&& 16267 \\\hline
F=29  & 1  & 6  & 30  & 77  & 241  & 337  & 850  & 1104  & {\bf 2249}  & {\bf 2744}  & {\bf 3200}  & {\bf 3584}  & {\bf 4096}  & {\bf 8192}  & {\bf 4096}  & {\bf 2048}  & {\bf 1024}  & {\bf 512}  & {\bf 256}  & {\bf 128}  & {\bf 64}  & {\bf 32}  & {\bf 16}  & {\bf 8}  & {\bf 4}  & {\bf 2}  & 1  & 0  & 1 && 34903 \\\hline
F=30  & 0  & 0  & 20  & 0  & 0  & 331  & 676  & 991  & {\bf 0}  & {\bf 2668}  & {\bf 3072}  & {\bf 3584}  & {\bf 4096}  & {\bf 0}  & {\bf 8192}  & {\bf 4096}  & {\bf 2048}  & {\bf 1024}  & {\bf 512}  & {\bf 256}  & {\bf 128}  & {\bf 64}  & {\bf 32}  & {\bf 16}  & {\bf 8}  & {\bf 4}  & {\bf 2}  & 1  & 0  & 1 & 31822 \\\hline
\hline\end{tabular}}}} 
  \caption{Number of semigroups with Frobenius number $F$ and multiplicity $m$, for $F$ from $1$ to $30$ and $2\leq m\leq F+1$.}
\label{taulaFm}
\end{table}

\section{Fibonacci-like growth with respect to the genus}
\label{s:ggrowth}

Kaplan proves in \cite{Kaplan}
that, if $N(m,g)$ is the number of numerical semigroups of multiplicity $m$ and genus $g$, then, under the hypothesis that $m$ and $g$ satisfy $2g \leq 3m-1$, it holds $$N (m - 1, g - 1) + N (m - 1, g - 2) = N (m, g).$$
Here we prove a more refined version of this result, by taking also the Frobenius number $F$ of a numerical semigroup into account.

The next lemma gives a well known formula for $n_{F,m,g}$ for the case in which $m>\frac{F}{2}.$

\begin{lemma}\label{fmgmgran}
  Let $F,m$ be nonnegative integers with $\frac{F}{2}<m\leq F-1$.
  Then for any integer $g$,
  $$n_{F,m,g}=\binom{F-m-1}{g-m}.$$
\end{lemma}

\begin{proof}
  If $m> \frac{F}{2}$, the semigroups with Frobenius number $F$, multiplicity $m$ are in bijection with the subsets of $[m+1,F-1]$, of size $g-m$, where the correspondence is given by assigning to a semigroup the subset of gaps between $m+1$ and $F-1$. This proves the formula.
  \end{proof}

\begin{theorem}\label{kaplanrefinat}
  Let $F,m$ be nonnegative integers such that
  $\frac{F+2}{3}\leq m\leq F-1$. It holds
$$n_{F,m,g}=n_{F-2,m-1,g-1}+n_{F-2,m-1,g-2}.$$
\end{theorem}

\begin{proof}
If $m>\frac{F}{2}$, 
by Lemma~\ref{fmgmgran},
\begin{eqnarray*}n_{F,m,g}&=&\binom{F-m-1}{g-m}\\&=&\binom{F-m-2}{g-m}+\binom{F-m-2}{g-m-1}\\&=&\binom{(F-2)-(m-1)-1}{(g-1)-(m-1)}+\binom{(F-2)-(m-1)-1}{(g-2)-(m-1)}
\\&=&n_{F-2,m-1,g-1}+n_{F-2,m-1,g-2}\end{eqnarray*}

Suppose now that $m\leq\frac{F}{2}$.
  Since $\frac{F+1}{3}\leq m\leq \frac{F}{2}$, 
    by Theorem~\ref{pot2},
  $$n_{F,m,g}=\sum_{A\in B_{d}}\binom{\nu_{d,m}(A)}{\gamma_{d,F}(A)-g}$$
    and, since
$\frac{(F-2)+1}{3}\leq (m-1)\leq \frac{(F-2)}{2}$,
    by Theorem~\ref{pot2},
  $$n_{F-2,m-1,g-1}=\sum_{A\in B_{d}}\binom{\nu_{d,m-1}(A)}{\gamma_{d,F-2}(A)-(g-1)}$$
    and
$$n_{F-2,m-1,g-2}=\sum_{A\in B_{d}}\binom{\nu_{d,m-1}(A)}{\gamma_{d,F-2}(A)-(g-2)}$$
    Now,
  \begin{eqnarray*}
    n_{F-2,m-1,g-1}+n_{F-2,m-1,g-2}&=&\sum_{A\in B_{d}}\binom{\nu_{d,m-1}(A)}{\gamma_{d,F-2}(A)-(g-1)}+\sum_{A\in B_{d}}\binom{\nu_{d,m-1}(A)}{\gamma_{d,F-2}(A)-(g-2)}\\
    &=&\sum_{A\in B_{d}}\left(\binom{\nu_{d,m-1}(A)}{\gamma_{d,F-2}(A)-(g-1)}+\binom{\nu_{d,m-1}(A)}{\gamma_{d,F-2}(A)-(g-2)}\right)\\
    &=&\sum_{A\in B_{d}}\binom{1+\nu_{d,m-1}(A)}{2+\gamma_{d,F-2}(A)-g}\\
    &=&\sum_{A\in B_{d}}\binom{1+(m-1)-2-\sigma_d(A)}{2+(F-2)-2-\#A-\sigma_d(A)-g}\\
    &=&\sum_{A\in B_{d}}\binom{m-2-\sigma_d(A)}{F-2-\#A-\sigma_d(A)-g}\\
    &=&\sum_{A\in B_{d}}\binom{\nu_{d,m}(A)}{\gamma_{d,F}(A)-g}=n_{F,m,g}
    \end{eqnarray*} 
\end{proof}

As an example,
$n_{24,9,17}+n_{24,9,18}=87+76=163=n_{26,10,19}$, which can be checked in Table~\ref{tF24} and Table~\ref{tF26}.

We notice that Kaplan's result is a consequence of Theorem~\ref{kaplanrefinat}. Indeed, suppose that $2g\leq3m-1$, then any numerical semigroup with multiplicity $m$ and genus $g$ will have Frobenius number $F\leq 2g-1\leq3m-2$, and, so, it will hold that $m\geq \frac{F+2}{3}$. Now, by Theorem~\ref{kaplanrefinat},
\begin{eqnarray*}N(m,g)&=&\sum_{F=g}^{2g-1}n_{F,m,g}=\sum_{F=g}^{2g-1}\left(n_{F-2,m-1,g-1}+n_{F-2,m-1,g-2}\right)\\
  &=&\sum_{F=g-2}^{2(g-1)-1}\left(n_{F,m-1,g-1}+n_{F,m-1,g-2}\right)\\
  &=&\sum_{F=g-2}^{2(g-1)-1}n_{F,m-1,g-1}+\sum_{F=g-2}^{2(g-2)-1+2}n_{F,m-1,g-2},\end{eqnarray*}
Now, since $n_{g-2,m-1,g-1}=0$, we have $\sum_{F=g-2}^{2(g-1)-1}n_{F,m-1,g-1}=\sum_{F=g-1}^{2(g-1)-1}n_{F,m-1,g-1}=N(m-1,g-1)$.
And, since $n_{2(g-2)-1+1,m-1,g-2}=n_{2(g-2)-1+2,m-1,g-2}=0$, we have
$\sum_{F=g-2}^{2(g-2)-1+2}n_{F,m-1,g-2}=\sum_{F=g-2}^{2(g-2)-1}n_{F,m-1,g-2}=N(m-1,g-2)$, and the result follows.

Notice that the hypotheses of the theorem are necessary. Indeed,
$n_{20,6,12}=7$, while $n_{18,5,11}+n_{18,5,10}=4+2.$
On the other hand, $n_{21,5,12}=5$, while $n_{19,4,11}+n_{19,4,10}=4+3$.
So, we have counterexamples with 
$n_{F,m,g}>n_{F-2,m-1,g-1}+n_{F-2,m-1,g-2}$ and counterexamples with $n_{F,m,g}<n_{F-2,m-1,g-1}+n_{F-2,m-1,g-2}.$

\begin{remark}\label{remkaplan}
  The second part of the previous proof, that is, the case when $F> 2m$, could have been proved in the language of Kunz coordinates, by means of the bijection used in Kaplan's proof \cite[Theorem 1]{Kaplan}.
  It is used there that semigroups with $\frac{F+2}{3}\leq m\leq F-1$ have Kunz coordinates $(k_1,\dots,k_{m-1})$ with $k_1,\dots,k_{m-2}\in\{1,2,3\}$ and $k_{m-1}\in \{1,2\}$. The genus of a semigroup with Kunz coordinates $(k_1,\dots,k_{m-1})$ is $k_1+\dots+k_{m-1}$.
  Let $K(m,g)$ denote the set of Kunz coordinates of semigroups with multiplicity $m$ and genus $g$.
  Now,
  if $\frac{F+2}{3}\leq m\leq F-1$, then the next mapping is a bijection.
  $$
  \begin{array}{rcl}
    K(m,g)&\to&K(m-1,g-1)\cup K(m-1,g-2)\\
(k_1,\dots,k_{m-1})&\mapsto&(k_1,\dots,k_{m-2})
  \end{array}$$
  This is due to the fact that if $(k_1,\dots,k_{m-1})\in K(m,g)$, with $k_1,\dots,k_{m-2}\in\{1,2,3\}$ and $k_{m-1}\in\{1,2\}$, then
  $(k_1,\dots,k_{m-2})$ is also a set of Kunz coordinates of length decreased by one and total weight decreased by $k_{m-2}$, which is either $1$ or $2$.
  On the other hand, if $(k_1,\dots,k_{m-2})\in K(m-1,g-1)$ with $k_1,\dots,k_{m-2}\in\{1,2,3\}$, then $(k_1,\dots,k_{m-2},1)\in K(m,g)$, while, if $(k_1,\dots,k_{m-2})\in K(m-1,g-2)$ with $k_1,\dots,k_{m-2}\in\{1,2,3\}$, then $(k_1,\dots,k_{m-2},2)\in K(m,g)$.
  Let us restrict now to the case of semigroups with $F>2m$. If a semigroup has Kunz coordinates $(k_1,\dots,k_{m-1})$, then the condition $F>2m$ implies that there must be at least one coordinate $k_i=3$. If $k_i$ is the coordinate with largest index $i$ such that $k_i=3$, then the Frobenius number of the corresponding semigroup is $2m+i$. Since $k_{m-1}\in\{1,2\}$, we deduce that $i<m-1$. The previous bijection will map $(k_1,\dots,k_{m-1})$ to $(k_1,\dots,k_{m-2})$ and the corresponding semigroups will have respectively Frobenius number $2m+i$ and $2(m-1)+i$, hence, decreasing by $2$.
  \end{remark}

\section*{Acknowledgement}

This research is supported by project PID2024-156636NB-C21 (MATSE), funded by MCIN/AEI/10.13039/501100011033/ FEDER, UE.

This work was partially supported by the Wallenberg AI, Autonomous Systems and
Software Program (WASP) funded by the Knut and Alice Wallenberg Foundation.
This research was conducted using the resources of High Performance Computing
Center North (HPC2N).
The computations were enabled by resources provided by the National Academic
Infrastructure for Supercomputing in Sweden (NAISS), partially funded by the
Swedish Research Council through grants agreement no. 2025-22-1579, NAISS
2026/3-170.

The authors would like to thank Nathan Kaplan for his interesting remarks, and specially for suggesting the alternative proof of Theorem~\ref{kaplanrefinat} in Remark~\ref{remkaplan}, by means of Kunz coordinates.



\end{document}